\documentclass[11pt,a4paper]{article}
\usepackage[T1]{fontenc}
\usepackage[utf8]{inputenc}
\usepackage{lmodern}
\usepackage{amsmath,amssymb,amsthm,mathtools}
\usepackage[margin=2.55cm]{geometry}
\usepackage{enumitem}
\usepackage{microtype}
\usepackage{authblk}
\usepackage[hidelinks]{hyperref}
\hypersetup{
  pdftitle={A SQUARE-ROOT COMPLEX INEQUALITY AND ITS INDUCED METRIC STRUCTURE},
  pdfauthor={Gangsong Leng and Lecheng Yang},
  pdfsubject={Integral complex inequalities and induced metric geometry},
  pdfkeywords={square-root complex inequality, circle-valued functions, metric geometry, Hausdorff dimension}
}
\numberwithin{equation}{section}

\newtheorem{theorem}{Theorem}[section]
\newtheorem{lemma}[theorem]{Lemma}

\newtheorem{proposition}[theorem]{Proposition}
\theoremstyle{remark}
\newtheorem{remark}[theorem]{Remark}

\newcommand{\T}{\mathbb T}
\newcommand{\C}{\mathbb C}

\title{A SQUARE-ROOT COMPLEX INEQUALITY AND ITS INDUCED METRIC STRUCTURE}
\author[1]{Gangsong Leng}
\author[2]{Lecheng Yang}
\affil[1]{School of Mathematical Sciences, East China Normal University, Shanghai, China\\
\texttt{lenggangsong@163.com}}
\affil[2]{Shanghai High School, Shanghai, China\\
\texttt{stanleyyoeng44@gmail.com}}

\date{}

\begin{document}
\maketitle

\begin{abstract}
Let $(\Omega,\mu)$ be a finite measure space with $M=\mu(\Omega)>0$. We investigate the integral form, stability, and metric geometry associated with a square-root complex. After proving the inequality and determining all equality cases, we analyze its phase stability near the intersection of the two branches of the equality set. In general phase directions, the quadratic term is precisely a Cauchy--Schwarz deficit; along the corresponding degenerate cone, the leading term is of fourth order and is strictly positive. A symmetric two-point example shows that the exponent four is unavoidable in any uniform distance-stability estimate. Finally, on the group of measurable circle-valued functions, we introduce the LY-metric
\[
 d_\mu(f,g)=\left|M-\int_\Omega f\overline g\,d\mu\right|^{1/2}.
\]
We prove that this metric is bi-invariant and complete, and that it induces the same topology as the $L^2$ metric. On finite-dimensional tori, we establish the optimality of the exponent $1/2$, derive explicit formulas for the intrinsic distance and geodesics, describe the anisotropic geometry and volume growth of small metric balls, and show that the Hausdorff dimension is $n+1$.
\end{abstract}



\section{Introduction}

Complex inequalities combine algebraic and geometric phenomena. Because the complex field has no order compatible with its field operations, many standard real-variable arguments are unavailable. At the same time, modulus, conjugation, and the multiplicative geometry of the closed unit disk provide effective substitutes.

Haoyu Wang proposed the following inequality in~\cite{Wang2026}. Let $n$ be a positive integer and let $a_1,\ldots,a_n,b_1,\ldots,b_n\in\C$ satisfy
\[
 |a_k|\le1,\qquad |b_k|\le1\qquad(1\le k\le n).
\]
Then
\begin{equation}\label{eq:Wang}
 2\left(\left|n-\sum_{k=1}^n a_k\right|+
 \left|n-\sum_{k=1}^n b_k\right|\right)
 \ge \left|n-\sum_{k=1}^n a_kb_k\right|.
\end{equation}
Thus the aggregate deviation of the product sequence from the constant sequence $1$ is controlled linearly by the corresponding deviations of the two factor sequences.

Building on this result, we obtained the square-root inequality
\begin{equation}\label{eq:CY-discrete}
\left|n-\sum_{k=1}^n a_kb_k\right|^{1/2}
\le
\left|n-\sum_{k=1}^n a_k\right|^{1/2}
+
\left|n-\sum_{k=1}^n b_k\right|^{1/2}.
\end{equation}
If
\[
 A=\left|n-\sum_{k=1}^n a_k\right|,
 \qquad
 B=\left|n-\sum_{k=1}^n b_k\right|,
\]
then $(\sqrt A+\sqrt B)^2\le2(A+B)$, so~\eqref{eq:CY-discrete} immediately implies~\eqref{eq:Wang}. Hence the square-root complex inequality strengthens Wang's inequality. More significantly, the exponent $1/2$ converts~\eqref{eq:CY-discrete} into a multiplicative triangle inequality and thereby gives rise to a natural metric on groups of circle-valued functions.

The remainder of the paper is organized as follows. Section~2 proves the integral square-root complex inequality and characterizes all equality cases. Section~3 studies phase stability near the intersection of the two branches of the equality set. We first identify, through lemmas, the quadratic leading term in general directions and its degenerate cone. The only main stability theorem then gives the exact fourth-order expansion and a strictly positive lower bound along the degenerate directions. A symmetric two-point example shows that the distance exponent cannot be lowered below four in a uniform geometric stability estimate. The deficit decomposition from the original proof is relegated to an auxiliary lemma. Section~4 studies the metric induced by the square-root complex inequality. We prove bi-invariance, completeness, and equivalence with the $L^2$ topology. On finite-dimensional tori we further establish the optimality of the exponent $1/2$, derive explicit formulas for the intrinsic distance and geodesics, analyze the anisotropy and small-ball volume growth, determine the Hausdorff dimension, and discuss a multiplicative parallelogram law.

\section{The integral form of the square-root complex inequality}

\begin{theorem}[Integral square-root complex inequality]\label{thm:integral-CY}
Let $(\Omega,\mu)$ be a finite measure space with $M=\mu(\Omega)>0$. If $f,g:\Omega\to\C$ are measurable and
\[
 |f|\le1,\qquad |g|\le1\quad\text{a.e.},
\]
then
\begin{equation}\label{eq:integral-CY}
\left|M-\int_\Omega fg\,d\mu\right|^{1/2}
\le
\left|M-\int_\Omega f\,d\mu\right|^{1/2}
+
\left|M-\int_\Omega g\,d\mu\right|^{1/2}.
\end{equation}
Equality holds if and only if $f=1$ a.e. or $g=1$ a.e.
\end{theorem}

\begin{proof}
Set
\[
X=\int_\Omega(1-f)\,d\mu,\qquad
Y=\int_\Omega(1-g)\,d\mu,
\]
and
\[
Z=\int_\Omega(1-f)(1-g)\,d\mu.
\]
The identity
\[
1-fg=(1-f)+(1-g)-(1-f)(1-g)
\]
gives
\[
M-\int_\Omega fg\,d\mu=X+Y-Z.
\]
Since $|f|\le1$, pointwise we have
\[
|1-f|^2=1+|f|^2-2\Re f\le2(1-\Re f),
\]
and hence
\[
\int_\Omega|1-f|^2\,d\mu\le2\Re X\le2|X|.
\]
Similarly, $\int|1-g|^2\,d\mu\le2|Y|$. By Cauchy--Schwarz,
\[
|Z|\le
\left(\int|1-f|^2\,d\mu\right)^{1/2}
\left(\int|1-g|^2\,d\mu\right)^{1/2}
\le2\sqrt{|X||Y|}.
\]
Therefore
\[
\left|M-\int_\Omega fg\,d\mu\right|
\le|X|+|Y|+2\sqrt{|X||Y|}
=(\sqrt{|X|}+\sqrt{|Y|})^2,
\]
and taking square roots yields~\eqref{eq:integral-CY}.

We now discuss equality. Sufficiency is immediate. Conversely, suppose equality holds. Then equality must hold at every step of the inequality chain used above. If $X=0$, then
\[
0=\Re X=\int_\Omega(1-\Re f)\,d\mu.
\]
The integrand is nonnegative, so $f=1$ a.e. The case $Y=0$ is analogous.

Assume now that $X,Y\ne0$. Equality forces $X$ and $Y$ to be positive real numbers, $|f|=|g|=1$ a.e., and equality in Cauchy--Schwarz. Thus there exists $\lambda\in\C$ such that
\[
1-f=\lambda(1-\overline g)\quad\text{a.e.}
\]
Equality in the final triangle inequality also implies that $Z$ and $X+Y$ point in opposite directions. Since
\[
Z=\lambda\int|1-g|^2\,d\mu,
\]
we conclude that $\lambda$ is a nonpositive real number. For circle-valued functions,
\[
\Re(1-f)=\frac12|1-f|^2,
\qquad
\Re(1-g)=\frac12|1-g|^2.
\]
If $\lambda<0$, the linear relation above gives
\[
\Re(1-f)=\lambda\Re(1-g)\le0,
\]
so $f=g=1$ a.e., contradicting $X,Y\ne0$. If $\lambda=0$, then $f=1$ a.e. immediately. Hence equality can occur only in the two cases stated in the theorem.
\end{proof}

\begin{remark}
If $\Omega=\{1,\ldots,n\}$ is equipped with counting measure, Theorem~\ref{thm:integral-CY} reduces to the finite-sum inequality~\eqref{eq:CY-discrete}.
\end{remark}

\section{Stability of the square-root complex inequality}

Throughout this section, $(\Omega,\mu)$ is a finite measure space with $M=\mu(\Omega)>0$. For measurable complex-valued functions satisfying $|f|,|g|\le1$ a.e., write
\[
A=\left|M-\int_\Omega f\,d\mu\right|,\qquad
B=\left|M-\int_\Omega g\,d\mu\right|,
\]
\[
T=\left|M-\int_\Omega fg\,d\mu\right|,
\qquad
\mathcal D(f,g)=(\sqrt A+\sqrt B)^2-T.
\]
Also set
\[
\mathbb B_\infty=\{h\in L^\infty(\Omega):|h|\le1\ \text{a.e.}\}
\]
and
\begin{equation}\label{eq:equality-set}
\mathcal E_\mu=
(\{1\}\times\mathbb B_\infty)
\cup
(\mathbb B_\infty\times\{1\}).
\end{equation}
By Theorem~\ref{thm:integral-CY}, $\mathcal E_\mu$ is exactly the equality set of the integral square-root complex inequality, and
\begin{equation}\label{eq:distance-equality-set}
\operatorname{dist}_{2}((f,g),\mathcal E_\mu)^2
=
\min\left\{
\int_\Omega|1-f|^2\,d\mu,
\int_\Omega|1-g|^2\,d\mu
\right\}.
\end{equation}

The principal result of this section concerns the most degenerate phase directions near the intersection of the two branches of the equality set. It identifies a strictly positive fourth-order term in those directions and shows that the exponent four is the critical one for any prospective uniform geometric stability estimate.

\begin{theorem}[Fourth-order stability along degenerate directions]\label{thm:local-distance-stability}
Let $(\Omega,\mu)$ be a probability space, and let $\varphi\in L^4(\Omega;\mathbb R)$ satisfy
\[
\int_\Omega\varphi\,d\mu=0,
\qquad \varphi\not\equiv0.
\]
Let $\lambda>0$, and define
\[
f_t=e^{it\varphi},\qquad g_t=e^{i\lambda t\varphi}.
\]
Write
\[
m_j=\int_\Omega\varphi^j\,d\mu\qquad(j=2,3,4).
\]
Then, as $t\to0$,
\begin{equation}\label{eq:quartic-expansion-main}
\mathcal D(f_t,g_t)
=
\lambda(1+\lambda)^2
\left(\frac{m_4}{8}-\frac{m_3^2}{12m_2}\right)t^4
+o(t^4).
\end{equation}
In particular,
\begin{equation}\label{eq:quartic-positive-main}
\mathcal D(f_t,g_t)
\ge
\frac{\lambda(1+\lambda)^2}{24}
\left(\int_\Omega\varphi^2\,d\mu\right)^2t^4
+o(t^4).
\end{equation}
\end{theorem}

Before proving the theorem, we record three lemmas. The first isolates the nonnegative deficits needed from the original proof. The second identifies the quadratic term for general pure-phase perturbations and characterizes its degenerate cone. The third gives the fourth-order expansion of the basic modulus functional.

\begin{lemma}[Basic deficit decomposition]\label{lem:basic-deficit}
Let
\[
X=\int_\Omega(1-f)\,d\mu,\qquad
Y=\int_\Omega(1-g)\,d\mu,
\]
\[
Z=\int_\Omega(1-f)(1-g)\,d\mu,
\]
and set
\[
P=\int_\Omega|1-f|^2\,d\mu,
\qquad
Q=\int_\Omega|1-g|^2\,d\mu.
\]
Then
\begin{align}\label{eq:basic-deficit-chain}
\mathcal D(f,g)
={}&\bigl(|X|+|Y|-|X+Y|\bigr)\notag\\
&+2\bigl(\sqrt{|X||Y|}-\sqrt{(\Re X)(\Re Y)}\bigr)\notag\\
&+\bigl(2\sqrt{(\Re X)(\Re Y)}-\sqrt{PQ}\bigr)\notag\\
&+\bigl(\sqrt{PQ}-|Z|\bigr)\notag\\
&+\bigl(|X+Y|+|Z|-|X+Y-Z|\bigr),
\end{align}
and all five terms on the right are nonnegative.
\end{lemma}

\begin{proof}
Since $|f|,|g|\le1$,
\[
P\le2\Re X,
\qquad Q\le2\Re Y,
\]
and Cauchy--Schwarz gives $|Z|\le\sqrt{PQ}$. Hence
\[
\begin{aligned}
|X+Y-Z|
&\le |X+Y|+|Z|\\
&\le |X+Y|+\sqrt{PQ}\\
&\le |X+Y|+2\sqrt{(\Re X)(\Re Y)}\\
&\le |X+Y|+2\sqrt{|X||Y|}\\
&\le |X|+|Y|+2\sqrt{|X||Y|}.
\end{aligned}
\]
Subtracting successive terms in this chain and summing the resulting deficits gives~\eqref{eq:basic-deficit-chain}.
\end{proof}

\begin{lemma}[Quadratic leading term]\label{lem:quadratic-main}
Let $(\Omega,\mu)$ be a probability space and let $\varphi,\psi\in L^4(\Omega;\mathbb R)$ satisfy
\[
\int_\Omega\varphi\,d\mu=
\int_\Omega\psi\,d\mu=0.
\]
Set $f_t=e^{it\varphi}$ and $g_t=e^{it\psi}$. Then, as $t\to0$,
\begin{equation}\label{eq:quadratic-expansion-main}
\mathcal D(f_t,g_t)
=t^2\bigl(\|\varphi\|_2\|\psi\|_2
-\langle\varphi,\psi\rangle\bigr)+o(t^2).
\end{equation}
The quadratic leading term is nonnegative, and it vanishes if and only if
\[
\psi=\lambda\varphi\quad\text{a.e.}
\]
for some $\lambda\ge0$.
\end{lemma}

\begin{proof}
By the zero-mean condition and Taylor expansion,
\[
\left|1-\int_\Omega e^{it\varphi}\,d\mu\right|
=\frac{t^2}{2}\|\varphi\|_2^2+o(t^2),
\]
and similarly for $\psi$. Moreover,
\[
\left|1-\int_\Omega e^{it(\varphi+\psi)}\,d\mu\right|
=\frac{t^2}{2}\|\varphi+\psi\|_2^2+o(t^2).
\]
Substitution into the definition of $\mathcal D$ gives~\eqref{eq:quadratic-expansion-main}. Nonnegativity and the equality condition follow from Cauchy--Schwarz; the proportionality factor must be nonnegative.
\end{proof}

\begin{lemma}[Fourth-order expansion of the modulus]\label{lem:quartic-main}
Let $(\Omega,\mu)$ be a probability space, and let $\varphi\in L^4(\Omega;\mathbb R)$ satisfy $\int_\Omega\varphi\,d\mu=0$. Write
\[
m_j=\int_\Omega\varphi^j\,d\mu\qquad(j=2,3,4).
\]
If $\varphi\not\equiv0$, then, as $s\to0$,
\begin{equation}\label{eq:modulus-fourth}
\left|1-\int_\Omega e^{is\varphi}\,d\mu\right|
=
\frac{m_2}{2}s^2+
\left(-\frac{m_4}{24}+\frac{m_3^2}{36m_2}\right)s^4
+o(s^4).
\end{equation}
\end{lemma}

\begin{proof}
The fourth-order Taylor expansion and $\varphi\in L^4$ give
\[
1-\int_\Omega e^{is\varphi}\,d\mu
=
\frac{m_2}{2}s^2+\frac{i m_3}{6}s^3
-\frac{m_4}{24}s^4+o(s^4).
\]
Since $m_2>0$, taking the modulus and expanding the square root yields~\eqref{eq:modulus-fourth}.
\end{proof}

\begin{proof}[Proof of Theorem~\ref{thm:local-distance-stability}]
Apply~\eqref{eq:modulus-fourth} with $s=t$, $s=\lambda t$, and $s=(1+\lambda)t$, and substitute the resulting expansions into the definition of $\mathcal D$. The quadratic terms cancel, and simplification gives~\eqref{eq:quartic-expansion-main}.

By Cauchy--Schwarz and monotonicity of $L^p$ norms on a probability space,
\[
m_3^2\le m_2m_4,
\qquad
m_4\ge m_2^2.
\]
Therefore
\[
\frac{m_4}{8}-\frac{m_3^2}{12m_2}
\ge\frac{m_4}{24}
\ge\frac{m_2^2}{24}.
\]
Substituting this into~\eqref{eq:quartic-expansion-main} yields~\eqref{eq:quartic-positive-main}.
\end{proof}

\begin{remark}[Necessity of the fourth-order distance]\label{rem:critical-order}
Lemma~\ref{lem:quadratic-main} shows that the positively proportional directions in Theorem~\ref{thm:local-distance-stability} are exactly the degenerate cone of the quadratic leading term. Conversely, consider a two-point probability space and let $f=g=e^{it}$ and $f=g=e^{-it}$ each occur with probability $1/2$. Then
\[
\mathcal D(f,g)=8\sin^4\frac t2,
\qquad
\int|1-f|^2=\int|1-g|^2=4\sin^2\frac t2.
\]
Hence, by~\eqref{eq:distance-equality-set},
\[
\mathcal D(f,g)
=\frac12\operatorname{dist}_{2}((f,g),\mathcal E_\mu)^4.
\]
Letting $t\to0$, we see that if a uniform estimate
\[
\mathcal D(f,g)\ge C\,\operatorname{dist}_{2}((f,g),\mathcal E_\mu)^p
\]
were to hold over all probability spaces and all admissible $f,g$, then necessarily $p\ge4$.
\end{remark}

\begin{remark}\label{rem:global-open}
Theorem~\ref{thm:local-distance-stability} completely describes the most degenerate phase model near the intersection of the two branches of the equality set: general directions are quadratically rigid, while the fourth-order leading term is strictly positive on the quadratic degenerate cone. To prove the global estimate
\[
\mathcal D(f,g)\ge\frac{c}{M}
\operatorname{dist}_{2}((f,g),\mathcal E_\mu)^4
\]
for all closed-unit-disk-valued functions, one would still need a uniform localization lemma that simultaneously handles radial deviations, nonzero mean phase, and nonuniform perturbations approaching $1$ only in the $L^2$ sense. We therefore do not state this global assertion as a theorem here.
\end{remark}

\section{The LY-metric}

This section shows that the square-root complex inequality not only gives an integral inequality, but also naturally defines a metric on a space of circle-valued functions. We first give the definition and basic properties on a general measure space, and then study intrinsic length, geodesics, metric balls, and Hausdorff dimension on finite-dimensional tori. For readability, all metric-geometric terms are explained when they first appear.

\subsection{Definition, topology, and completeness}

Let
\[
\T=\{z\in\C:|z|=1\}
\]
be the unit circle in the complex plane. We write $L^0(\Omega;\T)$ for the set of measurable circle-valued functions modulo equality almost everywhere. Under pointwise multiplication this set is an abelian group, with identity the constant function $1$ and inverse $\overline f$. Since $\mu(\Omega)<\infty$ and $|f|=1$ a.e., every $f\in L^0(\Omega;\T)$ belongs to $L^2(\Omega,\mu)$.

\begin{theorem}[The LY-metric]\label{thm:metric}
Let $(\Omega,\mu)$ be a finite measure space with $M=\mu(\Omega)>0$. For $f,g\in L^0(\Omega;\T)$ define
\begin{equation}\label{eq:metric}
 d_\mu(f,g)=\left|M-\int_\Omega f\overline g\,d\mu\right|^{1/2}.
\end{equation}
Then $d_\mu$ is a metric on $L^0(\Omega;\T)$. It is bi-invariant: for every $h\in L^0(\Omega;\T)$,
\begin{equation}\label{eq:bi-invariant}
 d_\mu(fh,gh)=d_\mu(hf,hg)=d_\mu(f,g).
\end{equation}
In particular, left and right multiplication are isometries.
\end{theorem}

\begin{proof}
Nonnegativity is immediate. Also,
\[
\int_\Omega g\overline f\,d\mu
=\overline{\int_\Omega f\overline g\,d\mu},
\]
so $d_\mu(f,g)=d_\mu(g,f)$.

If $d_\mu(f,g)=0$, then
\[
\int_\Omega f\overline g\,d\mu=M.
\]
Taking real parts gives
\[
\int_\Omega\bigl(1-\Re(f\overline g)\bigr)\,d\mu=0.
\]
Because $|f\overline g|=1$ a.e., the integrand is nonnegative. Thus $\Re(f\overline g)=1$ a.e., and consequently $f\overline g=1$ a.e.; hence $f=g$ a.e.

For the triangle inequality, set
\[
u=f\overline g,\qquad v=g\overline h.
\]
Then $|u|=|v|=1$ and $uv=f\overline h$. Theorem~\ref{thm:integral-CY} gives
\[
 d_\mu(f,h)\le d_\mu(f,g)+d_\mu(g,h).
\]
Finally, $|h|=1$ implies $(fh)\overline{(gh)}=f\overline g$, and~\eqref{eq:bi-invariant} follows.
\end{proof}

Define
\[
\Phi_\mu(f)=d_\mu(f,1)=\left|M-\int_\Omega f\,d\mu\right|^{1/2}.
\]
The function $\Phi_\mu$ is a length function based at the identity element of the group. Since
\[
d_\mu(f,g)=\Phi_\mu(f\overline g)
\]
and the square-root complex inequality yields
\[
\Phi_\mu(fg)\le\Phi_\mu(f)+\Phi_\mu(g),
\]
the square-root complex inequality is precisely the subadditivity of this multiplicative length function.

\begin{theorem}[Topological equivalence and completeness]\label{thm:L2}
For all $f,g\in L^0(\Omega;\T)$,
\begin{equation}\label{eq:l2-comparison}
\frac12\|f-g\|_{L^2(\mu)}^2
\le d_\mu(f,g)^2,
\qquad
d_\mu(f,g)^4\le M\|f-g\|_{L^2(\mu)}^2.
\end{equation}
Consequently, $d_\mu$ and the $L^2$ metric induce the same topology. Moreover, $(L^0(\Omega;\T),d_\mu)$ is complete.
\end{theorem}

\begin{proof}
Let $z=\int_\Omega f\overline g\,d\mu$. Then $|z|\le M$ and
\[
\|f-g\|_2^2=2(M-\Re z),
\qquad
d_\mu(f,g)^2=|M-z|.
\]
The first estimate follows from $M-\Re z\le|M-z|$. Since $|z|\le M$,
\[
|M-z|^2=M^2-2M\Re z+|z|^2
\le2M(M-\Re z),
\]
which gives the second estimate.

If $(f_j)$ is Cauchy with respect to $d_\mu$, then the first estimate shows that it is Cauchy in $L^2$. Let $f$ be its $L^2$ limit. Passing to a subsequence converging almost everywhere, we obtain $|f|=1$ a.e., so $f\in L^0(\Omega;\T)$. The second estimate then implies $d_\mu(f_j,f)\to0$. Thus the metric space is complete.
\end{proof}

\subsection{Finite-dimensional tori and the optimal exponent}

Let $\Omega=\{1,\ldots,n\}$ with counting measure. Then
\[
 \T^n=\underbrace{\T\times\cdots\times\T}_{n\ \text{factors}}
 =\{(z_1,\ldots,z_n)\in\C^n:|z_k|=1\}
\]
is called the $n$-dimensional torus. Each point can be written as $(e^{i\theta_1},\ldots,e^{i\theta_n})$, so $\T^n$ may also be viewed as the space of $n$ phase angles. Here the word ``torus'' refers to a Cartesian product of unit circles, not to a complex torus $\C^n/\Lambda$ in the sense of complex geometry.

On $\T^n$, the LY-metric becomes
\begin{equation}\label{eq:dn}
 d_n(u,v)=\left|n-\sum_{k=1}^n u_k\overline{v_k}\right|^{1/2}.
\end{equation}
It is invariant under coordinatewise multiplicative translations and under permutations of the coordinates.

\begin{theorem}[Optimality of the exponent $1/2$]\label{thm:optimal-exponent}
For $\alpha>0$, define
\[
 d_{n,\alpha}(u,v)=\left|n-\sum_{k=1}^n u_k\overline{v_k}\right|^\alpha.
\]
If $0<\alpha\le1/2$, then $d_{n,\alpha}$ is a metric for every positive integer $n$. If $\alpha>1/2$, then $d_{2,\alpha}$ already fails the triangle inequality. Hence the largest exponent valid uniformly on all finite-dimensional tori is $1/2$.
\end{theorem}

\begin{proof}
If $0<\alpha\le1/2$, then $d_{n,\alpha}=d_n^{2\alpha}$. If $d$ is a metric and $0<p\le1$, then $(a+b)^p\le a^p+b^p$ implies that $d^p$ is again a metric. This proves the first assertion.

Now let $\alpha>1/2$ and take
\[
x=(1,1),\qquad y=(e^{it},e^{-it}),\qquad z=(e^{2it},e^{-2it}).
\]
A direct computation gives
\[
d_{2,\alpha}(x,y)=d_{2,\alpha}(y,z)
=\bigl(2(1-\cos t)\bigr)^\alpha,
\]
whereas
\[
d_{2,\alpha}(x,z)=\bigl(2(1-\cos2t)\bigr)^\alpha.
\]
The triangle inequality would therefore imply
\[
\bigl(2(1+\cos t)\bigr)^\alpha\le2.
\]
As $t\to0$, the left-hand side tends to $4^\alpha>2$, so every sufficiently small $t>0$ yields a counterexample.
\end{proof}

\subsection{Intrinsic length and geodesics induced by \texorpdfstring{$d_n$}{dn}}

The original metric $d_n(u,v)$ gives the distance between two points directly. To study distances measured along curves in the torus, we introduce the path length induced by $d_n$.

For a continuous curve $\gamma:[0,1]\to\T^n$, define its $d_n$-length by
\[
L_{d_n}(\gamma)=\sup_{0=t_0<\cdots<t_N=1}
\sum_{j=1}^N d_n\bigl(\gamma(t_{j-1}),\gamma(t_j)\bigr),
\]
where the supremum is taken over all finite partitions of $[0,1]$. This is the standard metric length obtained from polygonal sums over finite partitions.

For $u,v\in\T^n$, define the intrinsic distance induced by $d_n$ by
\[
\rho_n(u,v)=\inf_\gamma L_{d_n}(\gamma),
\]
where the infimum is over all piecewise $C^1$ curves joining $u$ to $v$. If no joining curve has finite $d_n$-length, we set $\rho_n(u,v)=+\infty$. A curve that attains this infimum is called an intrinsic geodesic. Since the value may be $+\infty$, $\rho_n$ is, more precisely, a possibly extended intrinsic distance. The term ``intrinsic'' indicates that distance is obtained by minimizing lengths of paths in the space, rather than by evaluating the original two-point metric directly. On the total-phase leaves introduced below, this intrinsic distance agrees with the usual length distance of a flat Riemannian metric.

For $\omega\in\T$, define
\[
\mathcal L_\omega=
\left\{z=(z_1,\ldots,z_n)\in\T^n:\prod_{k=1}^n z_k=\omega\right\}.
\]
We call $\mathcal L_\omega$ a total-phase leaf. Indeed, if $z_k=e^{i\theta_k}$, then $\prod_kz_k=e^{i\sum_k\theta_k}$, so fixing $\omega$ is equivalent to fixing the total phase $\sum_k\theta_k$ modulo $2\pi$. The sets $\mathcal L_\omega$ form a partition of $\T^n$ into $(n-1)$-dimensional subtori.

\begin{theorem}[Intrinsic length and geodesics on the leaves]\label{thm:geodesics}
Every piecewise $C^1$ curve of finite $d_n$-length is contained entirely in one total-phase leaf $\mathcal L_\omega$.

Let $\gamma:[0,1]\to\mathcal L_\omega$ be piecewise $C^1$, and on each smooth subinterval choose continuous lifts of the phase functions so that
\[
\gamma(t)=\bigl(e^{i\theta_1(t)},\ldots,e^{i\theta_n(t)}\bigr).
\]
Then $\sum_{k=1}^n\theta_k'(t)=0$, and
\begin{equation}\label{eq:path-length-leaf}
L_{d_n}(\gamma)=
\frac1{\sqrt2}\int_0^1
\left(\sum_{k=1}^n\theta_k'(t)^2\right)^{1/2}dt.
\end{equation}

If $u,v\in\mathcal L_\omega$, choose real numbers $\delta_k$ such that $v_k=u_ke^{i\delta_k}$. Since $u$ and $v$ lie on the same leaf, $\sum_k\delta_k\in2\pi\mathbb Z$. Then
\begin{equation}\label{eq:intrinsic-distance}
\rho_n(u,v)=\frac1{\sqrt2}
\min_{\substack{m_k\in\mathbb Z\\
\sum_{k=1}^n(\delta_k+2\pi m_k)=0}}
\left(\sum_{k=1}^n(\delta_k+2\pi m_k)^2\right)^{1/2}.
\end{equation}
If the integers $m_k$ attain the minimum and
\[
\xi_k=\delta_k+2\pi m_k,
\]
then
\[
\gamma_k(t)=u_ke^{it\xi_k},\qquad 0\le t\le1,
\]
is an intrinsic geodesic from $u$ to $v$. Points on distinct leaves have intrinsic distance $+\infty$.
\end{theorem}

\begin{proof}
On a smooth subinterval, write
\[
\gamma(t)=\bigl(e^{i\theta_1(t)},\ldots,e^{i\theta_n(t)}\bigr)
\]
and put $\xi_k=\theta_k'(t)$. As $h\to0$,
\[
n-\sum_{k=1}^n\gamma_k(t+h)\overline{\gamma_k(t)}
=-ih\sum_{k=1}^n\xi_k+
\frac{h^2}{2}\sum_{k=1}^n\xi_k^2+o(h^2).
\]
If $\sum_k\xi_k\ne0$ at some point, continuity gives a neighborhood on which $|\sum_k\xi_k|$ stays bounded away from zero. Dividing that neighborhood into $N$ equal pieces, each increment has $d_n$-length of order $N^{-1/2}$, while the sum is of order $N^{1/2}$ and therefore diverges as $N\to\infty$. This contradicts finiteness of the curve length. Thus every finite-length curve satisfies
\[
\sum_{k=1}^n\theta_k'(t)=0
\]
on each smooth subinterval. Hence $\prod_k\gamma_k(t)$ is constant, and the curve lies in a single total-phase leaf.

On a leaf, $\sum_k\xi_k=0$, and therefore
\[
d_n\bigl(\gamma(t+h),\gamma(t)\bigr)
=|h|\left(\frac12\sum_{k=1}^n\xi_k^2\right)^{1/2}+o(|h|).
\]
The metric derivative is therefore the integrand in~\eqref{eq:path-length-leaf}, and the standard length formula for piecewise $C^1$ curves gives that identity.

After choosing a base point on a total-phase leaf, its universal covering space may be identified with the hyperplane
\[
H=\left\{x\in\mathbb R^n:\sum_{k=1}^n x_k=0\right\}.
\]
Formula~\eqref{eq:path-length-leaf} shows that intrinsic length on the leaf is $1/\sqrt2$ times Euclidean length on $H$. Hence the shortest curve between two endpoints is obtained from the shortest phase lift in $H$. All possible lifts are
\[
(\delta_1+2\pi m_1,\ldots,\delta_n+2\pi m_n),
\]
subject to the zero-sum condition. Choosing one with minimal Euclidean norm gives~\eqref{eq:intrinsic-distance}; projecting the corresponding straight segment in $H$ to the torus gives the stated geodesic.
\end{proof}

Consequently, each total-phase leaf, equipped with its intrinsic distance, is a flat $(n-1)$-dimensional torus. Here ``flat'' means that, on its universal covering space, lengths are computed using the ordinary Euclidean norm.

The original metric $d_n$ is generally not a length metric. A metric is called a length metric if the distance between any two points equals the infimum of lengths of curves joining them. Here one only has
\[
d_n(u,v)\le\rho_n(u,v),
\]
and the inequality may be strict. For example, if
\[
u=(1,1),\qquad v=(e^{it},e^{-it}),\qquad0<t<\pi,
\]
then
\[
d_2(u,v)=2\sin\frac t2<t=\rho_2(u,v).
\]
Thus $d_n$ is analogous to chord length on a circle, whereas $\rho_n$ is analogous to arc length measured along a leaf.

\subsection{Metric balls, volume growth, and Hausdorff dimension}
The ball in $\T^n$ under the LY-metric is called an LY-ball.
By multiplicative invariance, all LY-balls of the same radius are isometric, so it suffices to study balls centered at $\mathbf1=(1,\ldots,1)$. Let
\[
B_r(u)=\{v\in\T^n:d_n(u,v)<r\}
\]
denote the open metric ball, and define the centroid map
\[
\mathfrak m(z)=\frac1n\sum_{k=1}^nz_k.
\]
Then
\begin{equation}\label{eq:ball-centroid}
B_r(\mathbf1)=
\left\{z\in\T^n:
\mathfrak m(z)\in D\left(1,\frac{r^2}{n}\right)\right\},
\end{equation}
where $D(c,R)=\{w\in\C:|w-c|<R\}$ is the open disk in the complex plane with center $c$ and radius $R$. Thus an LY-ball is the inverse image of a planar disk under the centroid map.

Write $z_k=e^{i\theta_k}$ and decompose the phase vector orthogonally as
\[
\theta=a\frac{(1,\ldots,1)}{\sqrt n}+x,
\qquad x\perp(1,\ldots,1).
\]
Here $a$ is the component in the common-phase direction, while $x$ is the balanced component with zero total phase. As $\theta\to0$,
\begin{equation}\label{eq:anisotropic}
d_n(\mathbf1,e^{i\theta})^4
=na^2+\frac14|x|^4+o(a^2+|x|^4).
\end{equation}
Hence a small ball of radius $r$ has scale $r$ in the balanced directions but only scale $r^2$ in the common-phase direction. Thus the small-scale geometry is anisotropic.

\begin{theorem}[Small-ball volume and Hausdorff dimension]\label{thm:volume}
Let $m_n$ be normalized Haar measure on $\T^n$, namely the unique translation-invariant probability measure satisfying $m_n(\T^n)=1$. If $n\ge2$, then as $r\downarrow0$,
\begin{equation}\label{eq:volume-asymptotic}
m_n(B_r(u))=C_n r^{n+1}+o(r^{n+1}),
\end{equation}
where
\begin{equation}\label{eq:Cn}
C_n=
\frac{2^{(n-1)/2}\pi^{(n-1)/2}}
{(2\pi)^n\sqrt n\,\Gamma((n-1)/2)}
B\left(\frac{n-1}{4},\frac32\right).
\end{equation}
For $n=1$, one has $m_1(B_r)=\pi^{-1}r^2+o(r^2)$. Consequently,
\begin{equation}\label{eq:Hausdorff}
\dim_{\mathrm H}(\T^n,d_n)=n+1.
\end{equation}
Here $\dim_{\mathrm H}$ denotes Hausdorff dimension, which measures the effective dimension through the asymptotic behavior of coverings by small metric balls.
\end{theorem}

\begin{proof}
By~\eqref{eq:anisotropic}, the leading-order model of a small ball is
\[
na^2+\frac14|x|^4<r^4,
\qquad x\in\mathbb R^{n-1}.
\]
The Euclidean volume of this model ball is
\[
\frac2{\sqrt n}|\mathbb S^{n-2}|
\int_0^{\sqrt2r}\rho^{n-2}
\sqrt{r^4-\frac{\rho^4}{4}}\,d\rho.
\]
Substituting $\rho=\sqrt2rs$ and using
\[
\int_0^1s^{n-2}\sqrt{1-s^4}\,ds
=\frac14B\left(\frac{n-1}{4},\frac32\right)
\]
gives the constant $C_n$. Under the anisotropic scaling $a=r^2A$ and $x=rX$, the normalized remainder in~\eqref{eq:anisotropic} converges uniformly to zero on bounded sets. Inner and outer approximations by slightly dilated model balls then yield~\eqref{eq:volume-asymptotic}.

By multiplicative invariance, the asymptotic is independent of the center. Hence there exist $r_0,c,C>0$ such that for every $u\in\T^n$ and $0<r<r_0$,
\[
c r^{n+1}\le m_n(B_r(u))\le C r^{n+1}.
\]
This is small-scale $(n+1)$-Ahlfors regularity: every sufficiently small metric ball has measure comparable to the $(n+1)$st power of its radius. The standard mass-distribution and covering arguments then give the matching lower and upper bounds for the Hausdorff dimension, proving~\eqref{eq:Hausdorff}.
\end{proof}

\subsection{The parallelogram law}

Let 
\[
\ell_n(u)=d_n(u,\mathbf1)
\]
be the LY-length function.
In a linear space, a norm comes from an inner product if and only if it satisfies the parallelogram identity. Since $\T^n$ is a multiplicative group rather than a linear space, the natural analogue is to ask whether
\[
\ell_n(uv)^2+\ell_n(uv^{-1})^2
=2\ell_n(u)^2+2\ell_n(v)^2
\]
holds identically.

\begin{proposition}\label{prop:parallelogram}
The LY-length function does not in general satisfy this multiplicative parallelogram identity. Hence the global LY-length function is not Hilbertian in the sense of this quadratic identity.
\end{proposition}

\begin{proof}
For $n=1$, take $u=v=e^{it}$ with $0<t<\pi/2$. Then
\[
\ell_1(uv)^2+\ell_1(uv^{-1})^2
=|1-e^{2it}|=2\sin t,
\]
whereas
\[
2\ell_1(u)^2+2\ell_1(v)^2
=4|1-e^{it}|=8\sin\frac t2.
\]
Since $2\sin t<8\sin(t/2)$, the two sides are unequal.
\end{proof}

Nevertheless, on the balanced tangent space of each total-phase leaf,
\[
H=\left\{\xi\in\mathbb R^n:\sum_{k=1}^n\xi_k=0\right\},
\]
the intrinsic norm
\[
F(\xi)=\left(\frac12\sum_{k=1}^n\xi_k^2\right)^{1/2}
\]
is a constant multiple of the Euclidean norm and therefore satisfies the usual parallelogram identity. Thus the global geometry of the LY-metric is not Hilbertian, while its infinitesimal intrinsic geometry on each fixed-total-phase leaf is flat and Euclidean.

\medskip
In summary, the square-root complex inequality both strengthens Wang's complex inequality and carries a natural stability and metric-geometric structure. Near the intersection of the two branches of the equality set, the stability exhibits a quadratic--quartic stratification: general phase directions are quadratically rigid, whereas the fourth-order leading term is strictly positive on the quadratic degenerate cone. The symmetric two-point example shows that the distance exponent cannot be lowered below four in any uniform distance-stability estimate. The induced LY-metric has the same topology as the $L^2$ metric, but its small-ball volume grows like $r^{n+1}$. In the common-phase direction it exhibits square-root snowflake behavior, meaning that distance scales like the square root of the phase difference; within fixed-total-phase leaves, however, it reduces to flat Euclidean geometry. This coexistence of transverse snowflaking and leafwise Euclidean geometry is the most distinctive feature of the LY-metric.

\end{document}